\documentclass[11pt]{article}
\usepackage{amsmath}
\usepackage{amsthm}
\usepackage{pxfonts}
\usepackage{yfonts}
\usepackage{dsfont}
\usepackage{graphicx}
\usepackage{relsize}
\usepackage[hidelinks]{hyperref} 

\def\bel{\begin{equation}\label}
\def\eeq{\end{equation}}
\def\ds{\displaystyle}

\def\mt{\longrightarrow}
\def\v{\vskip 1em}

\def\Rec{{\bf R}}
\def\R{\mathbb R}

\def\C{\mathfrak{B}}

\def\N{{\bf N}}

\def\S{{\bf S}}

\def\E{{\bf E}}

\def\Q{{\bf Q}}

\def\A{{\bf A}}

\def\L{{\bf L}}
\def\U{{\bf U}}

\def\Hat{\widehat}

\def\vol{{\bf vol}}

\def\M{{\bf M}}

\def\Cup{{\bigcup}}

\def\alpha{\alphaup}
\def\beta{\betaup}
\def\gamma{\gammaup}
\def\delta{\deltaup}
\def\theta{\thetaup}
\def\xi{{\xiup}}
\def\eta{{\etaup}}
\def\tau{{\tauup}}
\def\rho{{\rhoup}}
\def\phi{{\phiup}}
\def\psi{{\psiup}}
\def\lambda{{\lambdaup}}
\def\omega{\omegaup}
\def\varphi{{\varphiup}}
\def\gamma{{\gammaup}}
\def\c{{\bf c}}

\def\w{{\bf w}}

\newtheorem{lemma}{Lemma}[section]

\begin{document}

\[\hbox{\LARGE{\bf Strong maximal function revisit on Heisenberg group}}\]

\[\hbox{Chuhan Sun}\]
\begin{abstract}
We prove the $\L^p$-boundedness of the strong maximal operator defined on a Heisenberg group $w.r.t$ an absolutely continuous measure satisfying the product $\A_\infty$-property. 
\end{abstract}

\section{Introduction}
\setcounter{equation}{0}
The study of certain operators  that  commute with a multi-parameter family of dilations,  dates back to the time of Jessen, Marcinkiewicz and Zygmund. A number of pioneering  results have been accomplished, for example
by 
C\'ordoba and Fefferman \cite{Cordoba-Fefferman},     Fefferman and Stein \cite{R.Fefferman-Stein}, Fefferman \cite{R.Fefferman},   M\"{u}ller, Ricci and Stein \cite{Muller-Ricci-Stein},
Journ\'{e} \cite{Journe'} and Pipher \cite{Pipher}.

In this paper, we consider the strong maximal function operator defined on a Heisenberg group with 
a multiplication law:
\bel{multiplication law}
\begin{array}{cc}\ds
(u,v,t)\odot(\xi,\eta,\tau)~=~\Big[u+\xi, v+\eta,t+\tau+\mu\big(u\cdot\eta-v\cdot\xi\big)\Big],\qquad\mu\in\R
\end{array}
\eeq
for every $(u,v,t)\in\R^n\times\R^n\times\R$ and $(\xi,\eta,\tau)^{-1}=(-\xi,-\eta,-\tau)\in\R^n\times\R^n\times\R$.

Denote $\Rec$ to be a rectangle in $\R^{2n+1}$ parallel to the coordinates. Moreover, 
\bel{R}
\Rec~=~\bigotimes_{i=1}^m\Q_i\times I~\subset~\bigotimes_{i=1}^m\R^{\N_i}\times\R,\qquad \hbox{\small{$\N_1+\N_2+\cdots+\N_m=2n$}}
\eeq
where
$\Q_i\subset\R^{\N_i},i=1,2,\ldots,m$ are cubes and  $I\subset\R$ is an interval.

A strong maximal operator $\M$ is initially defined on a Heisenberg group as 
\bel{M initial}
\begin{array}{lr}\ds
\M  f(u,v,t)~=~\sup_{\Rec\ni(0,0,0)}
\vol \{\Rec\}^{-1} \iiint_\Rec\left|f [(u,v,t)\odot(\xi,\eta,\tau)^{-1}]\right|d\xi d\eta d\tau.
\end{array}
\eeq
Let $\xi\mt u-\xi$, $\eta\mt v-\eta$ and $\tau\mt t-\tau$, $\M$ can be equivalently defined as
\bel{M equi}
\M f(u,v,t)~=~\sup_{\Rec\ni(u,v,t)}
\vol \{\Rec\}^{-1} \iiint_\Rec\left|f \left(\xi, \eta, \tau+\mu\left(u\cdot\eta-v\cdot\xi\right)\right)\right|d\xi d\eta d\tau.
\eeq
$\diamond$ {\small Throughout, $\C>0$ is regarded as a generic constant depending on its sub-indices.}

{\bf Theorem A: Christ, 1992}~~{\it Let $\M$ be defined in  (\ref{M initial}). We have
\bel{Result A}
\begin{array}{cc}\ds
\left\| \M f\right\|_{\L^q(\R^{2n+1})}~\leq~\C_{p}~\left\| f\right\|_{\L^p(\R^{2n+1})},\qquad 1<p<\infty.
\end{array}
\eeq}\\
The $\L^p$-boundedness of $\M$  defined on more general Nilpotent Lie groups can be found in the paper of Michael Christ \cite{Michael Christ}, in which the elegant argument is carried out using  a number of 'ingredients' developed previously by Ricci and Stein \cite{Ricci-Stein 1} and Christ \cite{Michael Christ 1}-\cite{Michael Christ 2}.

Our aim is to give a generalization of {\bf Theorem A} by defining $\M$ $w.r.t$ some appropriate absolutely continuous measure: $\omega(u,v)dudvdt$. Given any subset $E\subset\R^{2n+1}$, we write
\[\vol_\omega\{E\}~=~\iiint_E\omega(u,v)dudvdt.\]

Define
\bel{M}
\begin{array}{lr}\ds
\M_\omega f(u,v,t)~=~
\sup_{\Rec\ni(u,v,t)}\vol_\omega\{\Rec\}^{-1}\iiint_{\Rec}\left|f \left(\xi,\eta,\tau+\mu\left(u\cdot\xi-v\cdot\eta\right)\right)\right| \omega (\xi,\eta)d\xi d\eta d\tau.
\end{array}
\eeq
Let $(x_i,x_i^\dagger)\in\R^{\N_i}\times\R^{2n-\N_i}$ for every $i=1,2,\ldots,m$ where $\N_1+\N_2+\cdots+\N_m=2n$. 
We say 
\bel{A infinity}
\omega~\in~\bigotimes_{i=1}^m\A_\infty\left(\R^{2n}\right),
\eeq
if $\omega(\cdot, x_i^\dagger)$ satisfies the $\A_\infty$-property in $\R^{\N_i}$ for every $x_i^\dagger\in\R^{2n-\N_i}$.

$\diamond$ {\small For brevity, we abbreviate  
$\ds\left\|f\right\|_{\L^p(\R^{2n+1},\omega)}^p=\iiint_{\R^{2n+1}} |f(u,v,t)|^p\omega(u,v)dudvdt$, $1< p<\infty$}.
 \v
{\bf Theorem A*}~~{\it Let $\M_\omega$ be defined in (\ref{M}). Suppose $\omega\in\bigotimes_{i=1}^m\A_\infty\left(\R^{2n}\right)$. We have
\bel{Result One}
\left\| \M_\omega f\right\|_{\L^p\left(\R^{2n+1},\omega\right)}~\leq~\C_{p~\omega}~\left\| f\right\|_{\L^p\left(\R^{2n+1},\omega\right)},\qquad 1<p<\infty.
\eeq}

The proof of {\bf Theorem A*} is an application of 
 a multi-parameter covering lemma due to C\'{o}rdoba and Fefferman \cite{Cordoba-Fefferman}: Unlike the Vitali-type  covering lemma for cubes, there is no mutually disjointness between rectangles. Instead, this is replaced by the $\L^p$-norm of a summation of indicator functions supported on a  sequence of selected rectangles. We observe that  C\'{o}rdoba-Fefferman covering lemma is particularly useful to handle $\M_\omega$ defined on a Heisenberg group.  
 
 We prove {\bf Theorem A*} in Section 2.
 For the sake of self-containedness, we give a proof of the covering lemma within the required setting  in Section 3.

\section{Proof of Theorem A*}
\setcounter{equation}{0}
Consider  $(u, v, t)=(x_i,x_i^\dagger)\in\R^{\N_i}\times\R^{2n+1-\N_i}$ for every $i=1,2,\ldots,k$. We define
\bel{A infinity w}
\w\in\bigotimes_{i=1}^k\A_\infty\left(\R^{2n+1}\right), 
\eeq
where $\w(\cdot, x_i^\dagger)$ satisfies the $\A_\infty$-property in $\R^{\N_i}$ for every $x_i^\dagger\in\R^{2n+1-\N_i}$, $i=1,2,\cdots,k$.

For every subset $E\subset\R^{2n+1}$, we write
\bel{w meaure}
\vol_{\w}(E)~=~\iiint_{E}\w(u, v, t)dudvdt.
\eeq
Especially, if $\w=1$, then the measure becomes the Lebesgue measure, and
\bel{Lebesgue measure}
 \vol(E)=\vol_{{\bf 1}}(E)=\iiint_{E}{\bf 1} dudvdt.
 \eeq
Abbreviate
\[ \left\|f\right\|^p_{\L^p(\R^{2n+1}, \w)}=\iiint_{\R^{2n+1}}\left|f(u, v, t)\right|^p\w(u, v, t)dudvdt,\qquad1<p<\infty.\]

\v
{\bf C\'{o}rdoba-Fefferman covering lemma} ~~
{\it Let $\left\{\Rec_j\right\}_{j=1}^\infty$ be a collection of rectangles in $\R^{2n+1}$  parallel to the coordinates. Suppose $\w\in\bigotimes_{i=1}^k\A_\infty\left(\R^{2n+1}\right)$. There is a subsequence $\left\{\Hat{\Rec}_k\right\}_{k=1}^{\infty}$ such that
\bel{vol Compara 1}
\vol_{\w}\Bigg\{\Cup_j \Rec_j\Bigg\}~\lesssim~\vol_{\w}\Bigg\{ \Cup_k \Hat{\Rec}_k\Bigg\}
\eeq
and
\bel{indicator Sum 1}
\left\|\sum_k\chi_{\Hat{\Rec}_k}\right\|^p_{\L^p(\R^{2n+1},\w)}~\lesssim~\vol_{\w}\Bigg\{\Cup_k\Hat{\Rec}_k\Bigg\},\qquad \hbox{\small{$1<p<\infty$}}
\eeq
where $\chi$ is an indicator function.
}\\

We take $\w(u,v,t)=\omega(u,v)\cdot 1$, where $\omega\in\bigotimes_{i=1}^m\A_\infty\left(\R^{2n}\right)$. Then (\ref{vol Compara 1})-(\ref{indicator Sum 1}) become
\bel{vol Compara}
\vol_{\omega}\Bigg\{\Cup_j \Rec_j\Bigg\}~\lesssim~\vol_{\omega}\Bigg\{ \Cup_k \Hat{\Rec}_k\Bigg\}
\eeq
and
\bel{indicator Sum}
\left\|\sum_k\chi_{\Hat{\Rec}_k}\right\|^p_{\L^p(\R^{2n+1},\omega)}~\lesssim~\vol_{\omega}\Bigg\{\Cup_k\Hat{\Rec}_k\Bigg\},\qquad \hbox{\small{$1<p<\infty$}}.
\eeq

Let
\bel{U_lambda}
\U_\lambda~=~\Bigg\{ (x,y,t)\in\R^{n}\times\R^{n}\times\R\colon \M_\omega f(u,v,t)>\lambda\Bigg\}.
\eeq

Given any $(u,v,t)\in\U_\lambda$, there is a rectangle $\Rec_j\ni(u,v,t)$ such that
\bel{R_j est}
 \vol_{\omega}\{\Rec_j\}^{-1}\iiint_{\Rec_j}\left|f(\xi, \eta, \tau+\mu(u\cdot \eta-v\cdot \xi))\right| \omega (\xi,\eta)d\xi d\eta d\tau~>~{1\over 2} \lambda.
\eeq
Let $(u,v,t)$ run through the set $\U_\lambda$. We have
\[\U_\lambda\subset\Cup_j~\Rec_j.\]
By applying the covering lemma, we select a subsequence $\{\Hat{\Rec}_k\}_{k=1}^\infty$ from the union above and
\bel{union Rk size}
\begin{array}{lr}\ds
\vol_\omega\Bigg\{ \U_\lambda\Bigg\}~\lesssim~\vol_\omega\Bigg\{ \Cup_j \Rec_j\Bigg\}~\lesssim~ \vol_\omega\Bigg\{ \Cup_k \Hat{\Rec}_k\Bigg\}\qquad \hbox{\small{by (\ref{vol Compara}) }}
\\\\ \ds~~~~~~~~~~~~~~~~~
~\leq~\sum_k\vol_\omega \left\{ \Hat{\Rec}_k\right\} 
\\\\ \ds~~~~~~~~~~~~~~~~~
~\leq~\sum_k {2\over \lambda}\iiint_{\Hat{\Rec}_k}\left|f(\xi, \eta, \tau+\mu(u\cdot \eta-v\cdot \xi))\right|\omega (\xi,\eta) d\xi d\eta d\tau
\qquad
\hbox{\small{by (\ref{R_j est})}}.
\end{array}
\eeq
Furthermore, we find
\bel{union size}
\begin{array}{lr}\ds
\vol_\omega\Bigg\{ \Cup_k \Hat{\Rec}_k\Bigg\}~\lesssim~
\lambda^{-1}\sum_k \iiint_{\Hat{\Rec}_k}\left|f(\xi, \eta, \tau+\mu(u\cdot \eta-v\cdot \xi))\right| \omega (\xi,\eta) d\xi d\eta d\tau
\\\\ \ds~~~~~~~~~~~~~~~~~~~~~~~
~=~\lambda^{-1}\iiint_{\R^{2n+1}}\left|f(\xi, \eta, \tau+\mu(u\cdot \eta-v\cdot \xi))\sum_k\chi_{\Hat{\Rec}_k}(\xi,\eta,\tau)\right|\omega (\xi,\eta)d\xi d\eta d\tau
\\\\ \ds~~~~~~~~~~~~~~~~~~~~~~~
~\leq~\lambda^{-1} \left\{ \iiint_{\R^{2n+1}}\left|f(\xi, \eta, \tau+\mu(u\cdot \eta-v\cdot \xi))\right|^p \omega (\xi,\eta)d\xi d\eta d\tau\right\}^{1\over p}\left\|\sum_k\chi_{\Hat{\Rec}_k}\right\|_{\L^{p\over p-1}(\R^{2n+1},\omega)}
\\ \ds~~~~~~~~~~~~~~~~~~~~~~~~~~~~~~~~~~~~~~~~~~~~~~~~~~~~~~~~~~~~~~~~~~~~~~~~~~~~~~~~~~~~~~~~~~~~~~~~~~~~~~~~~~~~~~~~~~~~~~~~~~~~~
 \hbox{\small{ by H\"older inequality}}
\\\\ \ds~~~~~~~~~~~~~~~~~~~~~~~
~=~\lambda^{-1}  \left\{ \iint_{\R^{2n}}\left\| f(\xi, \eta, \cdot)\right\|_{\L^p(\R)}^p\omega(\xi,\eta) d\xi d\eta \right\}^{1\over p} \left\|\sum_k\chi_{\Hat{\Rec}_k}\right\|_{\L^{p\over p-1}(\R^{2n+1}, \omega)}
\\\\ \ds~~~~~~~~~~~~~~~~~~~~~~~
~\leq~\lambda^{-1}\left\| f\right\|_{\L^p(\R^{2n+1},\omega)}\vol_\omega\Bigg\{ \Cup_k \Hat{\Rec}_k\Bigg\}^{p-1\over p}
\qquad \hbox{\small{ by (\ref{indicator Sum})}}.
\end{array}
\eeq
This implies
\bel{union size'}
\vol_\omega\Bigg\{ \Cup_k \Hat{\Rec}_k\Bigg\}^{1\over p}~\lesssim~{1\over \lambda} \left\| f\right\|_{\L^p(\R^{2n+1},\omega)}.
\eeq
Let $\U_\lambda$ defined in (\ref{U_lambda}).
We have
\bel{weak-type gamma}
\begin{array}{lr}\ds
\vol_{\omega}\Bigg\{ (x,y,t)\in\R^n\times\R^n\times\R\colon \M_\omega f(u,v,t)>\lambda\Bigg\}^{1\over p} ~=~\vol_\omega\Bigg\{ \U_\lambda\Bigg\}^{1\over p}
\\\\ \ds~~~~~~~~~~~~~~~~~~~~~~~~~~~~~~~~~~~~~~~~~~~~~~~~~~~~~~~~~~~~~~~~~~~~~~~~~~~~~~~~~~~
~\lesssim~\vol_{\omega}\Bigg\{ \Cup_k \Hat{\Rec}_k\Bigg\}^{1\over p} \qquad\hbox{\small{by (\ref{union Rk size})}}
\\\\ \ds~~~~~~~~~~~~~~~~~~~~~~~~~~~~~~~~~~~~~~~~~~~~~~~~~~~~~~~~~~~~~~~~~~~~~~~~~~~~~~~~~~~
~\lesssim~{1\over \lambda} \left\|f\right\|_{\L^p(\R^{2n+1},\omega)}\qquad \hbox{\small{by (\ref{union size'}) }}.
\end{array}
\eeq
By using this weak type $(p,p)$-estimate and applying Marcinkiewicz interpolation theorem, we conclude that $\M_\omega$ is bounded on $\L^p(\R^{2n+1},\omega)$ for  $1<p<\infty$.

\section{Proof of the covering lemma}
\setcounter{equation}{0}

First, we need the following lemma proved by Fefferman \cite{R.Fefferman 1}:
\begin{lemma}\label{measure size}
If $\w\in\bigotimes_{i=1}^k\A_\infty\left(\R^{2n+1}\right)$, then $\w$ satisfies the following: If $\Rec\subset\R^{2n+1}$ is any rectangle with its sides parallel to the axes and $E\subset \Rec$ is such that $\vol\left(E\right)>\frac{1}{2}\vol\left(\Rec\right)$, then $\vol_{\w}(E)>\eta \vol_{\w}(\Rec)$, for some $\eta>0$.
\end{lemma}

{\bf Proof}. The proof is by induction on $k$. Assume the result is true for $k-1$. Consider a rectangle $\Rec$ as above, $\Rec=I\times J$ where $I$ is a rectangle in $\R^{\N_1}\times\R^{\N_2}\times\cdots\times\R^{\N_{k-1}}$, and $J$ is a cube in $\R^{\N_k}$, $\N_1+\N_2+\cdots+\N_k=2n+1$. 

Let $E\subset \Rec$ such that 
\bel{E/R}
\frac{\vol\left(E\right)}{\vol\left(\Rec\right)}>\frac{1}{2}.
\eeq
For each $x_k^\dagger\in I$, let $J_{x_k^\dagger}=\left\{(x_k^\dagger, x_k)~:~x_k\in J\right\}$. We claim that there exists $I'\subset I$ satisfying $\vol(I')\geq \varepsilon \vol(I)$, where $\varepsilon>0$ is small enough, such that for $x_k^\dagger\in I'$,
\bel{EcapJ}
\vol\left(E\cap J_{x_k^\dagger}\right)>\varepsilon \vol(J_{x_k^\dagger})=\varepsilon \vol(J).
\eeq
If not, we must have $\vol\left(\left\{x_k^\dagger\in I: \vol(E\cap J_{x_k^\dagger})>\varepsilon \vol(J)\right\}\right)\leq\varepsilon\vol\left(I\right)$. 

Divide $I$ into two parts: 
\bel{G B}
G=\left\{x_k^\dagger\in I: \vol (E\cap J_{x_k^\dagger})>\varepsilon \vol(J)\right\},\qquad B=\left\{x_k^\dagger \in I: \vol (E\cap J_{x_k^\dagger})\leq\varepsilon\vol(J)\right\}.
\eeq
Then
\bel{m(E)}
\begin{array}{lr}\ds
\vol(E)~=~\int_{B}\vol(E\cap J_{x_k^\dagger})dx_k^\dagger+\int_{G}\vol(E\cap J_{x_k^\dagger})dx_k^\dagger
\\\\ \ds~~~~~~~~~~~
~\leq~ \vol(B)\cdot\varepsilon\vol(J)+\vol(G)\cdot \vol(J)
\\\\ \ds~~~~~~~~~~~
~=~\Big[\vol(I)-\vol(G)\Big]\cdot\varepsilon\vol(J)+\vol(G)\cdot \vol(J)
\\\\ \ds~~~~~~~~~~~
~=~\varepsilon\vol(I)\cdot\vol(J)+(1-\varepsilon)\vol(G)\cdot \vol(J)
\end{array}
\eeq
Suppose $\vol(G)\leq\varepsilon\vol(I)$, then we further have
\bel{m E}
\begin{array}{lr}\ds
\vol(E)\leq\varepsilon\vol(I)\cdot \vol(J)+\left(1-\varepsilon\right)\cdot\varepsilon\vol(I)\cdot \vol(J)
\leq\left(2\varepsilon-\varepsilon^2\right)\vol(\Rec).
\end{array}
\eeq
We can choose $\varepsilon>0$ small enough such that $0<2\varepsilon-\varepsilon^2<\frac{1}{2}$, then
\bel{vol E}
\vol(E)<\frac{1}{2}\vol(\Rec),
\eeq
which is contradicted to $\frac{\vol(E)}{\vol(\Rec)}>\frac{1}{2}$.

Now, from (\ref{EcapJ}) and since $\w(x_k^\dagger,\cdot)$ is $\A^\infty$ in the $x_k$ variable for every $x_k^\dagger\in\R^{2n+1-\N_k}$, we have
\bel{int x_m}
\int_{E\cap J_{x_k^\dagger}}\w(x_k^\dagger, x_k)dx_k\geq\eta\int_{J_{x_k^\dagger}}\w(x_k^\dagger,x_k) dx_k, \qquad x_k^\dagger\in I'.
\eeq
But also if we fix any $x_k\in J$, then
\bel{int x_m'}
\int_{I'}\w(x_k^\dagger, x_k) dx_k^\dagger\geq\eta'\int_{I}\w(x_k^\dagger, x_k) dx_k^\dagger,
\eeq
by induction. It follows by integrating (\ref{int x_m}) in $x_k^\dagger\in I'$ that
\bel{int x 1}
\int_{E}\w(x) dx\geq \eta\int_{I'\times J}\w(x) dx,
\eeq
and integrating (\ref{int x_m'}) in $x_k\in J$ gives
\bel{int x 2}
\int_{I'\times J}\w(x) dx\geq\eta'\int_{I\times J}\w(x) dx=\int_{R}\w(x)dx.
\eeq
By combing (\ref{int x 1})-(\ref{int x 2}), we have
\bel{int x}
\int_{E}\w(x)dx\geq\eta\eta'\int_{R}\w(x)dx.
\eeq
Now we finish the proof of Lemma \ref{measure size}.

Then we continue to prove the covering lemma.
We re-arrange the order of $\{\Rec_j\}_{j=1}^\infty$ if necessary so that the cross-section volume of $\Rec_j$ in $\R^{\N_k}$ is decreasing as $j\mt\infty$. For brevity, we call it  $x_{k}$-cross section.
Denote $\Rec_j^*$ to be the rectangle co-centered with $\Rec_j$ having its $x_k$-cross section tripled and keeping the others same.
We select $\Hat{\Rec}_k$ from $\left\{\Rec_j\right\}_{j=1}^\infty$ as follows.

Let $\Hat{\Rec}_1=\Rec_1$. Having chosen $\Hat{\Rec}_1, \Hat{\Rec}_2,\ldots,\Hat{\Rec}_{N-1}$, we pick $\Hat{\Rec}_{N}$ as the first rectangle $\Rec$ on the list of $\Rec_j$'s after $\Hat{\Rec}_{N-1}$ so that
\bel{R condition}
\vol\left\{~\Rec\cap\left[\mathop{\bigcup}\limits_{k=1\atop{\Hat{\Rec}^*_k\cap \Rec\neq\emptyset}}^{N-1}\Hat{\Rec}^*_k\right]~\right\}~<~\frac{1}{2}\vol\left\{\Rec\right\}.
\eeq
Suppose $\Rec$ is an unselected rectangle. There is a positive number $M$ such that $\Rec$ is on the list of $\Rec_j$'s after $\Hat{\Rec}_M$ and
\bel{unselect}
\vol\left\{~\Rec\cap\left[\mathop{\bigcup}\limits_{k=1\atop{\Hat{\Rec}^*_k\cap \Rec\neq\emptyset}}^M \Hat{\Rec}^*_k\right]~\right\}~\ge~\frac{1}{2}\vol\left\{\Rec\right\}.
\eeq
Recall $\Hat{\Rec}^*_k$ whose $x_k$-cross section is tripled. Moreover, the $t$-side length of $\{\Rec_j\}_{j=1}^\infty$ is decreasing as $j\mt\infty$. On the $x_k$-cross section, the projection of $\Rec$  is covered by the projection of  the union inside (\ref{unselect}). 

Let $(x_1, x_2, \cdots, x_k)\in\Rec$. Then slice all rectangles with a plane through $(x_1, x_2, \cdots, x_k)$ perpendicular to the $x_k$-cross section. Denote $\S$, $\Hat{\S}_k$ and $\Hat{\S}^*_k$ to be the slices regarding to $\Rec$, $\Hat{\Rec}_k$ and $\Hat{\Rec}^*_k$. Consequently, (\ref{unselect}) implies
\bel{S condition >}
\vol\left\{~\S\cap\left[\mathop{\bigcup}\limits_{k=1}^M\Hat{\S}^*_k\right]~\right\}~\ge~\frac{1}{2}\vol\left\{\S\right\}.
\eeq

Since $\w\in\bigotimes_{i=1}^{k-1}\A_\infty\left(\R^{2n+1}\right)$, then $\w(\cdot, x_i^\dagger)$ satisfies the $\A_\infty$-property in $\R^{\N_i}$ for every $x_i^\dagger\in\R^{2n+1-\N_i}$, $i=1,2,\cdots,k$. This implies $\w(\cdot, x_k)\in\bigotimes_{i=1}^{k-1}\A_\infty\left(\R^{2n+1-\N_k}\right)$ for any fixed $x_k\in\R^{\N_k}$. 

For every subset $F\subset\R^{2n+1-\N_k}$, we define
\bel{F measure}
\vol_{\w(\cdot, x_k)}^{2n+1-\N_k}(F)=\int_{F}\w(x_k^\dagger, x_k)dx_k^\dagger, \qquad \textrm{for every } x_k\in\R^{\N_k}.
\eeq

Then by Lemma \ref{measure size}, we have 
\bel{S condition > lambda}
\vol_{\w(\cdot, x_k)}^{{2n+1-\N_k}}\left\{~\S\cap\left[\mathop{\bigcup}\limits_{k=1}^M\Hat{\S}^*_k\right]~\right\}~\ge~\eta\vol_{\w(\cdot, x_k)}^{{2n+1-\N_k}}\left\{\S\right\}, \qquad \eta>0.
\eeq 

Let $\M_{\w(\cdot,x_k)}^{2n+1-\N_k}$ be the strong maximal operator defined in $\R^{2n+1-\N_k}$. Observe that (\ref{S condition > lambda}) further implies
\bel{M>1/2}
\M_{\w(\cdot,x_k)}^{2n+1-\N_k} (\chi_{\Cup_k \Hat{\S}^*_k}) (x_k^\dagger)~>~\eta,\qquad x_k^\dagger\in\Cup_j  \S_j.
\eeq
It is well-known that $\M_{\w(\cdot, x_1^\dagger)}^{\N_1}$ is the Hardy-Littlewood maximal operator bounded on $\L^p(\R^{\N_1}, \w(\cdot, x_1^\dagger))$, for $1<p<\infty$. 
By induction, we may assume that $\M_{\w(\cdot, x_k)}^{2n+1-\N_k}$ is bounded on all $\L^p(\R^{2n+1-\N_k}, \w(\cdot,x_k))$, $1<p<\infty$.
From (\ref{S condition > lambda})-(\ref{M>1/2}), by applying the $\L^p$-boundedness of $\M_{\w(\cdot,x_k)}^{2n+1-\N_k}$, we find
\bel{bound}
\vol_{\w(\cdot, x_k)}^{2n+1-\N_k}\Bigg\{\Cup_j \S_j\Bigg\}~\lesssim~\vol_{\w(\cdot, x_k)}^{2n+1-\N_k}\Bigg\{ \Cup_k \Hat{\S}^*_k\Bigg\}.
\eeq
By using (\ref{bound}) and integrating in the $x_k$-coordinate, we have
\bel{vol compara}
\begin{array}{lr}\ds
\vol_{\w}\Bigg\{\Cup_j \hbox{\bf R}_j\Bigg\}~\lesssim~\vol_{\w}\Bigg\{ \Cup_k \Hat{\Rec}^*_k\Bigg\} ~\lesssim~\vol_{\w}\Bigg\{ \Cup_k \Hat{\Rec}_k\Bigg\}
\end{array}
\eeq
which is (\ref{vol Compara 1}).

On the other hand, (\ref{R condition}) implies
\bel{S condition <}
\vol\left\{~\Hat{\S}_{N}\cap\left[\mathop{\bigcup}\limits_{k=1}^{N-1}\Hat{\S}^*_k\right]~\right\}~<~\frac{1}{2}\vol\left\{\Hat{\S}_{N}\right\}.
\eeq
We are given that the measure $\w(\cdot, x_k)$ as above on the hyperplane ( if the hyperplane in $x_k=\c$, then $d\w(\cdot,\c)=\w(x_1, x_2, \cdots, x_{k-1}, \c)dx_1dx_2\cdots dx_{k-1}$ ) belongs to $\bigotimes_{i=1}^{k-1}\A_\infty\left(\R^{2n+1-\N_k}\right)$ for every  $x_k\in\R^{\N_k}$, so that for some $\eta>0$,
\bel{S condition < measure}
\vol_{\w(\cdot, x_k)}^{2n+1-\N_k}\left\{~\Hat{\S}_{N}\cap\left[\mathop{\bigcup}\limits_{k=1}^{N-1}\Hat{\S}^*_k\right]~\right\}~<~(1-\eta)\vol_{\w(\cdot, x_k)}^{2n+1-\N_k}\left\{\Hat{\S}_{N}\right\}.
\eeq
Denote $\Hat{\E}_N=\Hat{\S}_N\setminus\Cup_{k<N}\Hat{\S}_k$.  
From (\ref{S condition < measure}), we find 
\bel{vol E>S}
\vol_{\w(\cdot, x_k)}^{2n+1-\N_k}\left\{ \Hat{\E}_N\right\}~\geq~\eta\vol_{\w(\cdot, x_k)}^{2n+1-\N_k}\left\{\Hat{\S}_N\right\}.
\eeq
Let $\phi\in\L^{p\over p-1}(\R^{2n+1-\N_k},\w(\cdot, x_k))$ and $\left\|\phi\right\|_{\L^{p\over p-1}(\R^{2n+1-\N_k},\w(\cdot, x_k))}=1$. 
We have
\bel{integration}
\begin{array}{lr}\ds
\int_{\R^{2n+1-\N_k}}\phi\cdot\sum_k\chi_{\Hat{\S}_k} d\w(\cdot, x_k)~=~\sum_k\int_{\Hat{\S}_k}\phi d\w(\cdot, x_k)
\\\\ \ds~~~~~~~~~~~~~~~~~~~~~~~~~~~~~~~~~~~~~~~~~~~~~~~
~=~\sum_k \left\{\frac{1}{\vol_{\w(\cdot, x_k)}^{2n+1-\N_k}\{\Hat{\S}_k\}}\int_{\Hat{\S}_k}\phi d\w(\cdot, x_k)\right\} \vol_{\w(\cdot,x_k)}^{2n+1-\N_k}\left\{\Hat{\S}_k\right\}
\\\\ \ds~~~~~~~~~~~~~~~~~~~~~~~~~~~~~~~~~~~~~~~~~~~~~~~
~\leq~\sum_k \left\{\frac{1}{\vol_{\w(\cdot, x_k)}^{2n+1-\N_k}\{\Hat{\S}_k\}}\int_{\Hat{\S}_k}\phi d\w(\cdot, x_k) \right\} \cdot\frac{1}{\eta}\vol_{\w(\cdot, x_k)}^{2n+1-\N_k}\left\{\Hat{\E}_k\right\}\qquad\hbox{\small{by (\ref{vol E>S})}}
\\\\ \ds~~~~~~~~~~~~~~~~~~~~~~~~~~~~~~~~~~~~~~~~~~~~~~~
~\lesssim~\sum_k\int_{\Hat{\E}_k}\M_{\w(\cdot, x_k)}^{2n+1-\N_k}(\phi) d\w(\cdot, x_k)
\\\\ \ds~~~~~~~~~~~~~~~~~~~~~~~~~~~~~~~~~~~~~~~~~~~~~~~
~=~\int_{\Cup_k\Hat{\S}_k}\M_{\w(\cdot, x_k)}^{2n+1-\N_k}(\phi) d\w(\cdot, x_k).
\end{array}
\eeq
By applying H\"{o}lder inequality and the $\L^p$-boundedness of $\M_{\w(\cdot, x_k)}^{2n+1-\N_k}$, we find
\bel{boundedness of M}
\begin{array}{lr}\ds
\int_{\Cup_k\Hat{\S}_k}\M_{\w(\cdot, x_k)}^{2n+1-\N_k}(\phi) d\w(\cdot, x_k)~\leq~\left\|\M_{\w(\cdot, x_k)}^{2n+1-\N_k}(\phi)\right\|_{\L^{p\over p-1}(\R^{2n+1-\N_k}, \w(\cdot, x_k))}\vol_{\w(\cdot, x_k)}^{2n+1-\N_k}\Bigg\{\Cup_k\Hat{\S}_k\Bigg\}^{1\over p}
\\\\ \ds~~~~~~~~~~~~~~~~~~~~~~~~~~~~~~~~~~~~~~~~~~~~~~~
~\leq~\C_p~\vol_{\w(\cdot, x_k)}^{2n+1-\N_k}\Bigg\{\Cup_k\Hat{\S}_k\Bigg\}^{1\over p}.
\end{array}
\eeq
By substituting (\ref{boundedness of M}) to (\ref{integration}) and taking the supremum of $\phi$, we arrive at
\bel{p norm}
\left\|\sum_k\chi_{\Hat{\S}_k}\right\|_{\L^p(\R^{2n+1-\N_k},\w(\cdot, x_k))}~\leq~\C_p~\vol_{\w(\cdot, x_k)}^{2n+1-\N_k}
\Bigg\{\Cup_k\Hat{\S}_k\Bigg\}^{\frac{1}{p}}.
\eeq
Raising both sides of (\ref{p norm}) to the $p^{th}$ power and integrating over $x_k$ give us (\ref{indicator Sum 1}).

{\small School of Mathematical Sciences, Zhejiang University, Hangzhou, 310058, China}\\
{\small email: sunchuhan@zju.edu.cn}


\begin{thebibliography}{100}



\bibitem{Cordoba-Fefferman}{\small A.~C\'{o}rdoba and R.~Fefferman, {\it A geometric proof of the strong maximal theorem}, Annals of Mathematics {\bf 102}: 95-100, 1975.}



\bibitem{Michael Christ 1}{\small M.~Christ, {\it Hilbert transforms along curves. I. Nilpotent groups}, Annals of Mathematics {\bf 122}: no.3, 575-596, 1985.}

\bibitem{Michael Christ 2}{\small M.~Christ, {\it Hilbert transforms along curves, III. Rotational curvature}, preprint, 1984.}


\bibitem{Michael Christ}{\small M.~Christ, {\it The strong maximal function on a nilpotent group}, Transactions of the American Mathematical Society {\bf 331}: no.1, 1-13, 1992.}


\bibitem{R.Fefferman}{\small R.~Fefferman, {\it Harmonic Analysis on Product Spaces}, Annals of Mathematics {\bf 126}: no.1, 109-130, 1987.}

\bibitem{R.Fefferman 1}{\small R.~Fefferman, {\it Multiparameter Fourier Analysis} in Beijing Lectures in Harmonic Analysis, ed. by E. M. Stein, Princeton University Press, 1986.}




\bibitem{R.Fefferman-Stein}{\small R.~Fefferman and E.~M.~Stein, {\it Singular Integrals on Product Spaces},
Advances in Mathematics {\bf 45}: no.2, 117-143, 1982.}







\bibitem{Stein}{\small  E.~M.~Stein, {\it
Harmonic Analysis: Real-Variable Methods, Orthogonality and Oscillatory Integrals},
 Princeton University Press, 1993.}
 

\bibitem{Ricci-Stein 1}{\small F.~Ricci and E.~M.~Stein, {\it Harmonic analysis on nilpotent groups and singular integrals. II: Singular kernels
supported on submanifolds}, Journal of Functional Analysis {\bf 78}: 56-84, 1988.}








\bibitem{Muller-Ricci-Stein}{\small D. M\"{u}ller,~~F.~Ricci,~~E.~M.~Stein,  
{\it Marcinkiewicz Multipliers and Multi-parameter structures on Heisenberg (-type) group, I}, 
Inventiones Mathematicae {\bf 119}: no.2, 199-233, 1995.}



























\bibitem{Journe'}{\small J.~L.~Journ\'{e}, {\it Calder\'{o}n-Zygmund Operators on Product Spaces}, Revista Mathematica Iberoamericana {\bf 1}: no.3, 55-91, 1985. }

\bibitem{Pipher}{\small J.~Pipher, {\it Journ\'{e}'s Covering Lemma and Its Extension to Higher Dimensions}, Duke Mathematics Journal {\bf 53}: no.3, 683-690, 1986.}






\end{thebibliography}
\end{document}